\documentclass{article}
\usepackage{amsfonts,amssymb,mathrsfs,amscd,amsmath,amsthm}
\usepackage[russian, english]{babel}
\usepackage[T2A]{fontenc}
\usepackage[pdftex]{graphicx}
\usepackage{hyperref}
\usepackage{verbatim}
\usepackage{misccorr}
\usepackage{mathtools}
\usepackage{indentfirst}
\usepackage{alltt}
\usepackage{stmaryrd}
\usepackage{enumitem}
\usepackage{xurl}

\theoremstyle{thmstyleone}
\newtheorem{theorem}{Theorem}

\theoremstyle{thmstyletwo}%
\newtheorem{example}{Example}%

\theoremstyle{thmstylethree}%

\newcommand{\RomanNumeralCaps}[1]
{\MakeUppercase{\romannumeral #1}}

\begin{document}

%\title[Three-dimensional Lie algebras admitting Nijenhuis operators]{Three-dimensional Lie algebras admitting semisimple algebraic Nijenhuis operators}
\title{Three-dimensional Lie algebras admitting regular semisimple algebraic Nijenhuis operators}
\maketitle
\begin{center}
	\author{Zhikhareva Ekaterina Sergeevna\\Lomonosov Moscow State University}
\end{center}
\abstract{%
The aim of this paper is to classify all real and complex 3-dimensional Lie algebras admitting regular semisimple algebraic 
Nijenhuis operators. This problem is completely solved (see Theorems~\ref{t2} and \ref{t3}) by describing all Nijenhuis eigenbases 
for each 3-dimensional Lie algebra. It turns out that the answer is different in real and complex cases in the sence that there are 
real Lie algebras such that they do not admit an algebraic Nijenhuis operator, but their complexification admits such operators. 
An equally interesting question is to describe all algebraic Nijenhuis operators 
which are not equivalent by an automorphism of the Lie algebra. We give an answer to this question for some Lie algebras.
}



\section{Introduction}\label{sec1}

For arbitrary operator field $P$ the formula
$$
{\cal N}_P(\xi, \eta)=P[[P\xi,\eta]]+P[[\xi,P\eta]]-[[P\xi,P\eta]]-P^2[[\xi,\eta]]
$$
defines $(1, 2)$ tensor field, which is skew-symmetric in lower indices. Here $[[\,{,}\,]]$ stands for standard commutator 
of vector fields. This tensor field was introduced in \cite{nijenhuis} and is known as {\it Nijenhuis torsion}. 
We call $P$ a {\it Nijenhuis operator} if its Nijenhuis torsion identically vanishes.

Natural class of such operators are left-invariant Nijenhuis operators on a Lie group~$G$
(this means that an operator field is invariant under the left action of the group on itself). 
In this case the tensor field $P$ is completely defined by a linear operator $L: \mathfrak g \to \mathfrak g$, 
where $\mathfrak g$ is the Lie algebra of $G$. It satisfies the identity 
\begin{equation}\label{def}
    L[L\xi,\eta]+L[\xi,L\eta]-[L\xi,L\eta]-L^2[\xi,\eta] = 0,
\end{equation}
where $[\,{,}\,]$ is the Lie algebra commutator. In this case the operator is called an {\it algebraic Nijenhuis operator}. 
Such operators appear in the study of completely integrable systems (both classic and quantum) and in many questions in geometry 
and mathematical physics \cite{konyaev, sokolov, panasyuk, panasyuk2}.
Note that new methods for studying Nijenhuis operators and a program for their further research was recenly proposed in paper \cite{BMK}.

One may forget now about the group $G$ and define algebraic Nijenhuis operator by using the identity \eqref{def} only. 
We will follow this approach for the rest of the paper. Notice that if $\mathfrak g$ is Abelian, then any linear operator 
on it is an algebraic Nijenhuis operator.

\begin{example}
\rm{
Let $\mathfrak g = \mathfrak{gl}(n)$ and let $A, B$ be a pair of matrices with the property $A^2 = B^2 = \operatorname{Id}$. 
Then the operator $L$ given by the formula
$$
L(X) = AXB + BAX + BX - XB
$$
for arbitrary $X \in \mathfrak{gl}(n)$ is an algebraic Nijenhuis operator \cite{sokolov}. 
In this case one can show that for generic parameters $A,B$ the operator $L$ is diagonalizable, 
but its eigenvalues have algebraic multiplicity more than $1$ \cite{konyaev}. 
$\square$
}
\end{example}

\begin{example}
\rm{
Given a Lie algebra $\mathfrak g$, consider its one-dimensional central extension, that is $\mathfrak g \oplus \mathfrak c$, 
where $\mathfrak c$ is one-dimensional center. It is easy to show that 
for any $a \in \mathfrak g^*$ and $\eta \neq 0 \in \mathfrak c$ the formulas
$$
\begin{aligned}
& L \xi = a(\xi) \eta, \quad \xi \in \mathfrak g, \\    
& L \eta = 0
\end{aligned}
$$
define an algebraic Nijenhuis operator. 
Note that $L^2 = 0$, i.e., all the eigenvalues of this operator are zero. $\square$
}
\end{example}

A linear operator $L$ is called {\it regular semisimple} if all its eigenvalues are pairwise different and it is diagonalizable, 
that is it has an eigenbasis. The next theorem provides a characteristic property of Lie algebras admitting 
such Nijenhuis operators. It essentially can be derived from \cite{kosmann}, but to keep our work self-contained we provide it with proof (see Section~\ref{sec3}).

\begin{theorem}\label{t1} 
Real or complex Lie algebra $\mathfrak g$ admits a regular semisimple algebraic Nijenhuis operator if and only if there exists 
a basis $\eta_1, \dots, \eta_n$ such that
%$$
%[\eta_i, \eta_j] \in \langle \eta_i, \eta_j \rangle.
%[\eta_i, \eta_j] \in \operatorname{Span}\{\eta_i, \eta_j\}.
%$$
$$
[\eta_i, \eta_j] = \alpha\eta_i +\beta \eta_j, \quad  \alpha,\beta=\mathrm{const}.
$$
In other words, this condition means that any two basic vectors $\eta_i, \eta_j$ generate a two-dimensional 
Lie subalgebra in~$\mathfrak g$.
\end{theorem}

We call the basis from Theorem~\ref{t1} a {\it Nijenhuis eigenbasis}. In these terms, a Lie algebra admits 
a regular semisimple algebraic Nijenhuis operator if and only if it has at least one Nijenhuis eigenbasis. 

Notice that given a Nijenhuis eigenbasis one may construct many Nijenhuis operators taking different 
collections of eigenvalues. The common thing for these operators is that they pairwise commute. 
Thus, the set of Nijenhuis eigenbases is in one-to-one correspondence with the spaces of commuting algebraic Nijenhuis operators.

\begin{example}
\rm{
Theorem~\ref{t1} puts a strong restrictions on the corresponding Lie algebra. 
For example, let $\mathfrak g$ be a nilpotent Lie algebra admitting a Nijenhuis eigenbasis. Due to Theorem~\ref{t1}, 
any two vectors of this basis form a 2-dimensional subalgebra of $\mathfrak g$. At the same time, 
any subalgebra of a nilpotent algebra is itself nilpotent. Classification of 2-dimensional Lie algebras implies 
that if a 2-dimensional Lie algebra is nilpotent, then it is commutative. Thus, all elements of a Nijenhuis eigenbasis 
in $\mathfrak g$ commute. Therefore, any nilpotent Lie algebra $\mathfrak g$ admitting a regular semisimple algebraic Nijenhuis operator is commutative. 
$\square$
}   
\end{example}

\section{Results} \label{sec2}

Let us recall the well-known classification of 3-dimensional Lie algebras over $\mathbb R$ 
(originally due to Bianchi, we use notation from \cite{book}). We write only non-zero structural relations, 
the indices in our notation of the algebras are the dimension (in this paper we consider only 3-dimensional algebras) and the number in our list. 
The parameters $a$ for $\mathfrak a_{3, 7}$ and $b$ for $\mathfrak a_{3, 8}$ are arbitrary real numbers, 
which for different values yield non-isomorphic Lie algebras.
$$
\begin{aligned}
\begin{array}{lll}
 & \RomanNumeralCaps{1}& \mathfrak a_{3, 1} \quad Abelian,\\
 & \RomanNumeralCaps{2}& \mathfrak a_{3, 2} \quad [\eta_2, \eta_3] = \eta_1,\\
 & \RomanNumeralCaps{9}& \mathfrak a_{3, 3} \quad [\eta_1, \eta_2] = \eta_3, \quad [\eta_2, \eta_3] = \eta_1, \quad [\eta_3, \eta_1] = \eta_2,\\
 & \RomanNumeralCaps{8}& \mathfrak a_{3, 4} \quad  [\eta_1, \eta_2] = - \eta_3,  \quad  [\eta_2, \eta_3] = \eta_1, \quad  [\eta_3, \eta_1] = \eta_2,\\
 & \RomanNumeralCaps{5}& \mathfrak a_{3, 5} \quad  [\eta_1, \eta_2] = \eta_2, \quad [\eta_3, \eta_1]= - \eta_3, \\
 & \RomanNumeralCaps{4}& \mathfrak a_{3, 6} \quad  [\eta_1, \eta_2] = \eta_2 + \eta_3, \quad [\eta_3, \eta_1] = - \eta_3, \\
 & \RomanNumeralCaps{7}& \mathfrak a_{3, 7} \quad [\eta_1, \eta_2] = a \eta_2 + \eta_3, \quad [\eta_3, \eta_1] = \eta_2 - a \eta_3, \\
 & \RomanNumeralCaps{6}& \mathfrak a_{3, 8}\quad [\eta_1,\eta_2] = b \eta_2 - \eta_3, \quad [\eta_3,  \eta_1] = \eta_2 - b\eta_3. \\
\end{array}
\end{aligned}
$$

In particular, in the above list the Lie algebras $\mathfrak a_{3, 3}$ and $\mathfrak a_{3, 4}$ are isomorphic 
to 3-dimensional semisimple Lie algebras ${\rm so}(3,\mathbb R)$ and ${\rm sl}(2,\mathbb R)$ respectively.

The next theorem lists all the real Lie algebras which satisfy the property from Theorem~\ref{t1}.

\begin{theorem}\label{t2}
The only real non-Abelian three-dimensional Lie algebras that admit Nijenhuis eigenbases are 
$\mathfrak a_{3, 4}, \mathfrak a_{3, 5}$, and $\mathfrak a_{3, 8}$. The examples of Nijenhuis eigenbases $\zeta_1,\zeta_2,\zeta_3$
for each of them are given below (the vectors $\zeta_i$ are written in terms of the bases in the classification above).
\begin{itemize}
\item 
In $\mathfrak a_{3, 4}$ an example of a Nijenhuis eigenbasis $\zeta_1, \zeta_2, \zeta_3$ is given by 
$$
\zeta_1 = - \eta_1 + \eta_2 + \eta_3, \quad \zeta_2 = \eta_1 + \eta_2 - \eta_3, \quad \zeta_3 = 2\eta_1.
$$
\item 
In $\mathfrak a_{3, 5}$ any basis is a Nijenhuis eigenbasis
\item 
In $\mathfrak a_{3, 8}$ an example of a Nijenhuis eigenbasis $\zeta_1, \zeta_2, \zeta_3$ is given by 
$$
\zeta_1 = \eta_1 + \eta_3, \quad \zeta_2 = \eta_1 + \eta_2, \quad \zeta_3 = \eta_2 + \eta_3.
$$
\end{itemize}
\end{theorem}

We consider two Nijenhuis eigenbases equivalent if they are related by rescaling of vectors, 
that is $\zeta_i \to \lambda_i \zeta_i$ for some collection of non-zero constants $\lambda_i$. 
Further, when we say basis, we would actually mean the equivalence class of bases.

In case of 3-dimensional real Lie algebras the moduli space of Nijenhuis eigenbases can be described in a purely 
geometric manner (see the proof of Theorem \ref{t2} in Section~\ref{sec4}).

\begin{itemize}
    \item For algebra $\mathfrak a_{3, 4}$ the set of Nijenjhuis eigenbases is in one-to-one correspondence 
with triples of non-collinear points of $\mathbb RP^2$ lying on the given oval;
    \item For algebra $\mathfrak a_{3, 5}$ the set of Nijenhuis eigenbases is in one-to-one correspondence 
with triples of non-collinear points of $\mathbb RP^2$;
    \item For algebra $\mathfrak a_{3, 8}$ the set of Nijenhuis eigenbases is in one-to-one correspondence 
with triples of non-collinear points of $\mathbb RP^2$ lying on two given intersecting lines.
\end{itemize}

The classification of complex 3-dimensional Lie algebras is similar to the classification of real 3-dimensional 
Lie algebras, except that the Lie algebras $\mathfrak a_{3, 3}$ and $\mathfrak a_{3, 4}$ become isomorphic, 
and the Lie algebras $\mathfrak a_{3, 7}$ and $\mathfrak a_{3, 8}$ become parts of a single family of Lie algebras.

The following is the list of the 3-dimensional Lie algebras over $\mathbb C$ 
$$
\begin{aligned}
	\begin{array}{lll}
	& \RomanNumeralCaps{1}& \mathfrak b_{3, 1} \quad Abelian,\\
	& \RomanNumeralCaps{2}& \mathfrak b_{3, 2} \quad [\eta_2, \eta_3] = \eta_1,\\
	& \RomanNumeralCaps{9} \simeq \RomanNumeralCaps{8} & \mathfrak b_{3, 3} \quad [\eta_1, \eta_2] = \eta_3, \quad [\eta_2, \eta_3] = \eta_1, \quad [\eta_3, \eta_1] = \eta_2,\\
	&\RomanNumeralCaps{5}& \mathfrak b_{3, 4} \quad  [\eta_1, \eta_2] = \eta_2, \quad [\eta_3, \eta_1]= - \eta_3 \\
	&\RomanNumeralCaps{4}& \mathfrak b_{3, 5} \quad  [\eta_1, \eta_2] = \eta_2 + \eta_3, \quad [\eta_3, \eta_1] = - \eta_3, \\
	&\RomanNumeralCaps{7},\RomanNumeralCaps{6} & \mathfrak b_{3, 6} \quad [\eta_1, \eta_2] = a \eta_2 + \eta_3, \quad [\eta_3, \eta_1] = \eta_2 - a \eta_3, \qquad a \in \mathbb{C}\\
	\end{array}
\end{aligned}
$$

\begin{theorem}\label{t3}
The only complex non-Abelian three-dimensional Lie algebras that admit Nijenhuis eigenbases are $\mathfrak b_{3, 3}, \mathfrak b_{3, 4}$, and $\mathfrak b_{3, 6}$. 
\end{theorem}

The examples of Nijenhuis eigenbases for each complex Lie algebra are written in the last part of the proof of Theorem \ref{t3} (see Section~\ref{sec5}). 
%The main observation is that the lists in case of $\mathbb C$ and $\mathbb R$ differ.

One can note that lists in Theorems~\ref{t2} and \ref{t3} (in real and complex cases respectively) are similar, but they have differences 
in the following sence: for example, ${\rm so}(3,\mathbb R)$ does not admit a Nijenhuis eigenbasis, but ${\rm so}(3,\mathbb C)$
has such bases. Here both in real and complex cases this (orthogonal) Lie algebras are denoted in the above lists by $\RomanNumeralCaps{9}$.

%%%%%%%%%%%%%%%%%%%%%%%%%%%%%%%%%%%%%%%%%%%%%%%%%%%%%%%%%%%%%%%%%

\section{Proof of Theorem \ref{t1}} \label{sec3}

First, let us show that the existence of a Nijenhuis eigenbasis is necessary.
Consider a regular semisimple algebraic Nijenhuis operator $L$ on $\mathfrak g$ and prove that its eigenbasis
is a Nijenhuis eigenbasis. Let $\eta_i, \eta_j$ be a pair of eigenvectors with eigenvalues $\lambda_i, \lambda_j$.  
From condition (1) we get
$$
\begin{aligned}
0 & = L[L\eta_i,\eta_j]+L[\eta_i,L\eta_j]-[L\eta_i,L\eta_j]-L^2[\eta_i,\eta_j] = \\
& = \lambda_i L[\eta_i, \eta_j] + \lambda_j L [\eta_i, \eta_j] - \lambda_i \lambda_j [\eta_i, \eta_j] - L^2 [\eta_i, \eta_j] 
= - (L {-} \lambda_i \operatorname{Id})(L {-} \lambda_j \operatorname{Id}) [\eta_i, \eta_j].
\end{aligned}
$$
Since the operator $L$ is regular semisimple, the kernel $\operatorname{Ker} \big(L - \lambda_i \operatorname{Id}\big)$ is spanned 
by $\eta_i$ and the operators $(L-\lambda_i \operatorname{Id})$ and $(L-\lambda_j \operatorname{Id})$ pairwise commute. 
Thus, $[\eta_i, \eta_j]$ lies in the two-dimensional space spanned by $\eta_i, \eta_j$. 
So, we see that the existence of a Nijenhuis eigenbasis is necessary.
	
Now assume that at least one Nijenhuis eigenbasis exists. We pick arbitrary collection of pairwise distinct 
numbers $\lambda_1, \dots, \lambda_n$ and define 
$$
L \eta_i = \lambda_i \eta_i.
$$
It is enough to check the identity \eqref{def} in any basis. Performing the same calculations in the inverse direction, we get 
that \eqref{def} holds and the existence of a Nijenhuis eigenbasis is sufficient as well.

The theorem is proved.

%%%%%%%%%%%

\section{Proof of Theorem \ref{t2}}  \label{sec4}

Fixing a basis $\eta_1, \eta_2, \eta_3$ in a 3-dimensional Lie algebra $\mathfrak g$, we have
$$
[\eta_i, \eta_j] = c^k_{ij} \eta_k
$$
with structure constants $c^k_{ij}$. Consider a pair of vectors $\alpha,\beta\in\mathfrak g$:
$$
\begin{aligned}
\alpha &= a^1 \eta_1 + a^2 \eta_2 + a^3 \eta_3, \\
\beta  &= b^1 \eta_1 + b^2 \eta_2 + b^3 \eta_3. \\
\end{aligned}
$$
They form a subalgebra if and only if $\alpha \wedge \beta \wedge [\alpha, \beta] = 0$. Introducing
$$
M_1 = a^2\cdot b^3 - a^3\cdot b^2, \quad M_2 = a^3\cdot b^1 - a^1\cdot b^3, \quad M_3 = a^1\cdot b^2 - a^2 \cdot b^1,
$$
we get that this identity is equivalent to
\begin{equation}\label{xt1}
\begin{aligned}
0  = M^2_1 \cdot c^1_{23}+ M^2_2 \cdot c^2_{31}+M^2_3 \cdot c^3_{12} & + M_1M_2 \cdot (c^1_{31}+c^2_{23}) + \\
& + M_2M_3 \cdot (c^2_{12}+c^3_{31})+ M_1M_3 \cdot (c^1_{12}+c^3_{23}).
\end{aligned}    
\end{equation}
Regarding $\mathfrak g$ as a vector space equipped with standard Euclidean bilinear form in the basis $\eta_i$,
we can treat $M = (M_1, M_2, M_3)$ as the vector orthogonal to the two-dimensional space spanned by $\alpha, \beta$. 
In addition, it satisfies homogeneous quadratic equation \eqref{xt1}.

In these terms, the existence of a Nijenhuis eigenbasis $\zeta_1, \zeta_2, \zeta_3$ is equivalent 
to the existence of three linearly independent solutions for equation \eqref{xt1}. 
Substituting the structure constants for each 3-dimensional Lie algebra from the classification, we get the following forms 
of \eqref{xt1}
\begin{equation*}
\begin{array}{lll}
& \mathfrak a_{3, 2} \quad M_1^2 = 0\\ 
& \mathfrak a_{3, 3} \quad M_1^2+M_2^2+M_3^2=0\\ 
& \mathfrak a_{3, 4} \quad M_1^2+M_2^2-M_3^2=0 \\
& \mathfrak a_{3, 5} \quad \text{none}\\
& \mathfrak a_{3, 6} \quad M_3^2=0\\
& \mathfrak a_{3, 7} \quad M_2^2+M_3^2=0\\
& \mathfrak a_{3, 8} \quad M_2^2-M_3^2=0.\\
\end{array}
\end{equation*}
We see that the only algebras for which there are three linearly independent solutions of \eqref{xt1}
are $\mathfrak a_{3, 4}, \mathfrak a_{3, 5}$, and $\mathfrak a_{3, 8}$. For each of this algebras, except for $\mathfrak a_{3, 5}$, 
we provide an example of a Nijenhuis eigenbasis.

\begin{itemize}
\item 
For algebra $\mathfrak a_{3, 4}$ consider the basis
$$
\zeta_1 = - \eta_1 + \eta_2 + \eta_3, \quad \zeta_2 = \eta_1 + \eta_2 - \eta_3, \quad \zeta_3 = 2\eta_1.
$$
The structure constants in this basis are
$$
[\zeta_1, \zeta_2] = \zeta_1 - \zeta_2, \quad [\zeta_2, \zeta_3] = \zeta_3 - 2 \zeta_2,\quad [\zeta_3, \zeta_1] = -2 \zeta_1 - \zeta_3.
$$
\item 
For algebra $\mathfrak a_{3, 8}$ consider the basis
$$
\zeta_1 = \eta_1 + \eta_3, \quad \zeta_2 = \eta_1 + \eta_2, \quad \zeta_3 = \eta_2 + \eta_3.
$$
The structure constants in this basis are
$$
[\zeta_1, \zeta_2] = (1 + b) (\zeta_2 - \zeta_1),\quad [\zeta_2, \zeta_3] = (b - 1) \zeta_3, \quad [\zeta_3, \zeta_1] = (1 - b) \zeta_3.
$$
\end{itemize}

Thus, the theorem is proved.

%%%%%%%%%

\section{Proof of Theorem \ref{t3}}  \label{sec5}

The proof in the complex case is similar to the proof in the real case.

Fixing a basis $\eta_1, \eta_2, \eta_3$ in a complex 3-dimensional Lie algebra $\mathfrak g$ and
considering $\alpha,\beta\in\mathfrak g$ and $M = (M_1, M_2, M_3)$, we obtain exactly the same equation~\eqref{xt1} 
as in the real case. Then, substituting the structural constants for each complex 3-dimensional Lie algebra 
from the classification, we get the following forms of \eqref{xt1}
\begin{equation*}
\begin{array}{lll}
& \mathfrak b_{3, 2} \quad M_1^2 = 0,\\ 
& \mathfrak b_{3, 3} \quad M_1^2+M_2^2+M_3^2=0,\\ 
& \mathfrak b_{3, 4} \quad \text{none}\\
& \mathfrak b_{3, 5} \quad M_3^2=0\\
& \mathfrak b_{3, 6} \quad M_2^2+M_3^2=0\\
\end{array}
\end{equation*}
The only difference from the real case here is that the equations $M_1^2+M_2^2+M_3^2=0$ and $M_2^2+M_3^2=0$
now have three linearly independent solutions. So, we see that the only complex Lie algebras
that have Nijenhuis eigenbases are $\mathfrak b_{3, 3}$, $\mathfrak b_{3, 4}$, and $\mathfrak b_{3, 6}$. 
For each of them, except for $\mathfrak b_{3, 4}$, we provide an example of a Nijenhuis eigenbasis.
\begin{itemize}
\item 
For algebra $\mathfrak b_{3, 3}$ consider the basis
$$
\zeta_1 = - \eta_1 + \eta_2 + \eta_3, \quad \zeta_2 = \eta_1 + \eta_2 - \eta_3, \quad \zeta_3 = 2\eta_1.
$$
The structure constants in this basis are
$$
[\zeta_1, \zeta_2] = \zeta_1 - \zeta_2, \quad [\zeta_2, \zeta_3] = \zeta_3 - 2 \zeta_2,\quad [\zeta_3, \zeta_1] = -2 \zeta_1 - \zeta_3.
$$
\item 
For algebra $\mathfrak b_{3, 6}$ consider the basis
$$
\zeta_1 = \eta_1 + \eta_3, \quad \zeta_2 = \eta_1 + \eta_2, \quad \zeta_3 = \eta_2 + \eta_3.
$$
The structure constants in this basis are
$$
[\zeta_1, \zeta_2] = (1 + a) (\zeta_2 - \zeta_1),\quad [\zeta_2, \zeta_3] = (a - 1) \zeta_3, \quad [\zeta_3, \zeta_1] = (1 - a) \zeta_3.
$$
\end{itemize}

Thus, the theorem is proved.

\section{Remark on equivalence of algebraic Nijenhuis operators}

Theorems~\ref{t2} and~\ref{t3} give the classification of all real and complex 3-dimensional Lie algebras admitting regular
semisimple algebraic Nijenhuis operators and the description of all Nijenhuis eigenbases related to these operators.

In the case when a Lie algebra $\mathfrak g$ admits algebraic Nijenhuis operators, it is also natural to consider them up to
an automorphism of $\mathfrak g$. In other words, it is interesting to investigate for each Lie algebra how many 
really ``different'' regular semisimple algebraic Nijenhuis operators there are on it. Of course, each such operator has
a collection of eigenvalues, but from the geometric point of view they are not important. So, we can say that
two regular semisimple algebraic Nijenhuis operators $L$  and $L'$ on a 3-dimensional Lie algebra $\mathfrak g$ are equivalent 
if they are defined (in the sense of Theorem~\ref{t1}) by Nijenhuis eigenbases $\eta_1,\eta_2,\eta_3$ and $\eta'_1,\eta'_2,\eta'_3$ 
such that $\eta'_i=\mu_i\Phi(\eta_i)$ for some constants $\mu_i$ and an automorphism $\Phi:\mathfrak g\to\mathfrak g$.

The classification of algebraic Nijenhuis operators on 3-dimensional Lie algebras up to this equivalence
will be presented by the author in a further paper. Now let us only mentioned here some examples concerning this question.

1) It is evident that in the case of Abelian Lie algebra all bases are Nijenhuis bases and all of them are equivalent.

2) It turns out that for the Lie algebra ${\rm sl(2)}$ 
(real or complex, i.e., $\mathfrak a_{3, 4}$ or $\mathfrak b_{3, 3}$ in our notation), all Nijenhuis eigenbases
are also equivalent. This follows from the fact, which can be obtained by direct calculations, that in an arbitrary Nijenhuis eigenbasis $\zeta_1,\zeta_2,\zeta_3$ 
for the Lie algebra ${\rm sl(2)}$ the commutators have the form 
$$
[\zeta_1,\zeta_2]=B\zeta_1+A\zeta_2, \quad [\zeta_2,\zeta_3]=C\zeta_2+B\zeta_3, \quad [\zeta_3,\zeta_1]=A\zeta_3+C\zeta_1,
$$
where $A,B,C$ are arbitrary nonzero constants. 
Therefore, multiplying vectors $\zeta_1,\zeta_2,\zeta_3$ by the constants $\frac{A'}A,\frac{B'}B,\frac{C'}C$ respectively one can 
obtain an arbitrary triple of constants $A',B',C'$ from the given structural constants $A,B,C$. 

3) As follows from Theorems~\ref{t2} and~\ref{t3}, for the Lie algebra of type $\RomanNumeralCaps{5}$ 
(i.e., $\mathfrak a_{3, 5}$ or $\mathfrak b_{3, 4}$ in our notation) all bases are Nijenhuis bases, but
not all of them can be taken one to another by an automorphism and multiplication by constants. This can be seen, 
for example, from the fact that in this Lie algebra the commutant does not coinside with the algebra itself.
Therefore, in this case, unlike the previous examples, a regular semisimple algebraic Nijenhuis operator is not unique.

\smallskip

The author expresses gratitude to her scientific supervisor A.\,A.~Oshemkov and to A.\,Yu.~Konyaev
for usefull discussions under work on this paper, and also to A.\,V.~Bolsinov for valuable remarks.


\begin{thebibliography}{9}
\bibitem{nijenhuis} 
A. Nijenhuis, $X_n$-forming sets of eigenvectors, Nederl. Akad. Wetensch. Proc. Ser.~A 54: Indag. Math. 13 (1951).\\ 
https://doi.org/10.1016/S1385-7258(51)50028-8
\bibitem{konyaev} 
A. Yu. Konyaev, Completeness of commutative Sokolov--Odesskii subalgebras and Nijenhuis operators on gl(n), Sb. Math., 211:4 (2020), 583--593.\\ 
https://doi.org/10.1070/sm9251
\bibitem{sokolov} 
A. V. Odesskii, V. V. Sokolov, Integrable matrix equations related to pairs of compatible associative algebras, J. Phys. A 39:40 (2006), 12447--12456.\\ 
https://doi.org/10.1088/0305-4470/39/40/011
\bibitem{panasyuk} 
A. Panasyuk, Algebraic Nijenhuis operators and Kronecker Poisson pencils, Differential Geom. Appl. 24:5 (2006), 482--491.\\ 
https://doi.org/10.1016/j.difgeo.2006.05.009
\bibitem{panasyuk2} 
A. Panasyuk, A. Szereszewski, Webs, Nijenhuis operators, and heavenly PDEs, Classical and Quantum Gravity, 40:23 (2023), 235003.\\ 
https://doi.org/10.48550/arXiv.2211.03197
\bibitem{BMK}
A. Bolsinov, V. Matveev, A. Konyaev, Nijenhuis Geometry, Advances in Mathematics, 394 (2022), 108001.\\
https://doi.org/10.1016/j.aim.2021.108001 
\bibitem{kosmann} 
Y. Kosmann-Schwarzbach, F. Magri, Poisson--Nijenhuis structures. Ann. Inst. Henri Poincare, 53 (1990), 35--81.
\bibitem{book} 
B. A. Dubrovin, A. T. Fomenko, S. P. Novikov, 
Modern Geometry --- Methods and Applications. Part I. The Geometry of Surfaces, Transformation Groups, and Fields, Springer, New York (1984).
\end{thebibliography}
\end{document}